\documentclass{article}
\usepackage[utf8]{inputenc}

\usepackage{amsmath}
\usepackage{amssymb}
\usepackage{amsthm}
\usepackage{amsfonts}
\usepackage{mathtools}
\usepackage{blkarray}
\usepackage{tikz}
\usepackage{tikz-cd}

\usepackage{enumitem}
\usepackage[toc]{appendix}

\usepackage[authordate]{biblatex-chicago}

\usepackage[anythingbreaks]{breakurl} 
\usepackage{hyperref}
\usepackage{cleveref}
\crefname{section}{§\hspace{-0.1cm}}{§§}
\Crefname{section}{§}{§§}

\title{Cartesian Frames}
\author{Scott Garrabrant \\ Daniel A. Herrmann \\ Josiah Lopez-Wild}
\date{September 2021}

\begin{document}
\maketitle

\begin{abstract}
    We introduce a novel framework, the theory of Cartesian frames (CF), that gives powerful tools for manipulating sets of acts. The CF framework takes as its most fundamental building block that an agent can freely choose from a set of available actions. The framework uses the mathematics of Chu spaces to develop a calculus of those sets of actions, how those actions change at various levels of description, and how different agents' actions can combine when agents work in concert. We discuss how this framework might provide an illuminating perspective on issues in decision theory and formal epistemology. 
    
\end{abstract}

\section{Background and Motivation}
\label{back}

Acts are basic in decision theory. Once we have a set of acts available to an agent, we can then use decision theory to evaluate the choice-worthiness of the different acts. 

Much of the formal work philosophers have done to understand practical rationality focuses on this decision-theoretic \textit{evaluation} of acts (see for example \cite{jeffrey1965logic};  \cite{nozick1969newcomb}; \cite{savage1972foundations}; \cite{skyrms1984pragmatics}; \cite{skyrms1990dynamics}; \cite{joyce1999foundations}; \cite{ahmed2014evidence}). In this paper we take a step back, and provide a formal framework---the theory of Cartesian frames---for thinking about the acts \textit{themselves} in a more fine-grained way. 

The building blocks of the framework are \textit{Cartesian frames}.\footnote{This paper presents a fragment of the CF framework. The full description of the framework, originally developed by Garrabrant (\citeyear{garrabrant2020cf}), appeared on the AI Alignment Forum here: \url{https://www.alignmentforum.org/s/2A7rrZ4ySx6R8mfoT}.} Each Cartesian frame is a \textit{frame} in the sense that it is a particular way of carving up the world into an \textit{agent} and an \textit{environment}. Each way of carving up the world leads to a different frame. The frames are \textit{Cartesian} in the sense that, in any given frame, the agent is formally separate from the environment, and we think of the agent of that frame as able to freely choose between multiple acts. This fits well with Richard Jeffrey's standard definition that an ``act is then a proposition which is within the agent's power to make true if he pleases" (\citeyear{jeffrey1965logic}, p. 84).

In this paper we focus on two ways in which the CF framework extends our ability to pose and answer questions about agency. The first is that it gives a rigorous and precise way to think about acts at different levels of description. The second is that it gives us a systematic way to think about how different agents' acts can combine to give us more powerful agents (think of a group of people coming together to form a society), and how to decompose agents into subagents (think of starting with a society, and then realizing we can think of it as a collection of agents working together). 

Though the focus of this paper is the framework itself, we provide here a brief motivation for the framework.
\\
\\
\noindent \textit{Decisions at different levels of description.} The theory of Cartesian frames provides a powerful set of tools for expressing decision problems at different levels of description and for characterizing the relationships between the sets of acts available at those different levels of description. Understanding decision problems at different levels of granularity has already born fruit in decision theory. Consider a problem with which decision theorists are familiar: Savage's \textit{problem of small worlds} (\cite{savage1972foundations}). When thinking about how to represent a decision problem, Savage makes the distinction between \textit{small} and \textit{large} worlds (Chapter 5, Section 5, \citeyear{savage1972foundations}). A small world decision problem contains states, acts, and outcomes that are known to the agent, whereas a large world problem does not. Intuitively, the small world is one level of description, a coarser level that the agent can cognitively manage, whereas the large world problem is a more fine grained level, perhaps one which the agent lacks the conceptual resources to consider (see Joyce (\citeyear{joyce1999foundations}), p. 73 for a careful statement of problem). The problem of small worlds is the following: what guarantee can we have that the best act in a small world would also be the best in a larger world? One approach to this problem is to consider moving between different levels of description of a decision problem. Joyce takes such an approach, focusing on the process of refining decision problems (\cite{joyce1999foundations}, pp. 73-77). Given the expressive power of the theory of Cartesian frames, we expect that it will pay philosophical rent by allowing us to flexibly manipulate levels of description of decision problems, as well as connect these levels of description up to notions of subagency, to which we now turn. 
\\
\\

\noindent \textit{Subagency.} A novel contribution of the CF framework is its ability to model how different agents' actions are interrelated. This brings to the foreground the notion of a \textit{subagent}. Consider a corporation. At one level of description it is useful to think of the corporation as a single agent taking actions in the world.\footnote{For a philosophical discussions of this, see List, Pettit, et al. (\citeyear{list2011group}).} At another level of description, however, we can identify different agents that are parts of the bigger agent: for example, a sales department, a finance department, and so on. Yet another level down we can identify different agents: a personal assistant, a janitor, a CEO, a sales rep, etc. In this case we would say that the sales department is a \textit{subagent} of the corporation, and (for example) the personal assistant might be a subagent of both the sales department and the corporation, but not of the finance department. The theory of Cartesian frames has the resources to formally model these relationships. Level of description plays a certain role here, but it is not always from \textit{one} agent's perspective; at one level of description there might be one agent in the choice situation, at another there might be more than one.  
 
One area where we believe this rich calculus for subagency will be valuable is in exploring the connection between individual and group rationality, which is familiar to philosophers. The CF framework gives natural ways to consider the options agents have when working together, and indeed when it makes sense to consider them subagents of one composite agent (such as a group). Thus, insofar as this question relies on individuals coming together to form a group, the CF framework can help us by providing a precise and detailed way to think of this process. 

Consider, for example, the independence thesis, which is that (i) rational agents might form irrational groups and (ii) irrational agents might form rational groups.\footnote{For a thorough discussion of this thesis and how it relates to other positions see Mayo-Wilson, Zollman, and Danks (\citeyear{independencethesis}).} One way philosophers have tried to evaluate this thesis is to look at whether learning strategies that are successful in isolation are successful in group settings.\footnote{Mayo-Wilson, Zollman, and Danks show that successful individual strategies can be unsuccessful in group settings (\citeyear{independencethesis}, \citeyear{mayo2013wisdom}).} In this case a group is considered rational if all of its members are successful; the group itself does not have any actions available to it. Cartesian frames give a different way to approach this question that looks at the \textit{group} as an agent, since it gives operations that allow us to combine different agents into one \textit{superagent} composed of those agents. Thus the same standard of rationality (for example, whether or not an agent maximizes expected utility) can be applied at both the group  and the individual level. 

Note that this CF approach also differs from how similar questions are investigated in social choice theory. For example, Arrow's theorem asserts that there are groups of rational individuals with certain preferences over a set of alternatives that cannot be aggregated to a rational group preference over that same set of alternatives in a way that satisfies certain desirable conditions.\footnote{For an excellent presentation and discussion of Arrow's result see chapter 14 of Luce and Raiffa's \textit{Games and Decisions} (\citeyear{luce1989games}).} In this context we want to know whether or not we can construct a group preference that meets some conditions. The CF approach from the previous paragraph is different from this in two ways. The first is that it operates at the level of \textit{actions} instead of preferences. The second is that the set of options available to the group (the actions of the superagent) can be different from the actions available to the subagents, whereas in the social choice context the individual agents' preferences and the group's preferences are all defined over the same domain. Thus, the theory of Cartesian frames offers a novel way to investigate the relationship between group and individual rationality. 

Furthermore, such a formal framework allows us to ask new questions: are the values of the subagents (the individuals) of this agent (the group) aligned with the values of the agent? How sensitive to misaligned subagents is a superagent's ability to take successful action? Which subagents are most important to keep aligned? Such questions mirror important questions in the field of artificial intelligence. See Gabriel (\citeyear{gabriel2020artificial}) for a philosophical discussion of alignment in the context of AI. 

Precise notions of subagency might also give us better traction on understanding preexisting agent-based models used to tackle certain philosophical questions. For example, one thriving area of investigation in the philosophy of language is the evolution of compositionality---how is it that we humans developed such a rich, compositional language? Lewis and Skyrms took important first steps to answering this question (and more basic ones about simpler languages) in their respective works \textit{Convention} (\citeyear{lewis2008convention}) and \textit{Signals} (\citeyear{skyrms2010signals}). More recently a number of authors have used \textit{hierarchical models} to investigate this question (\cite{barrett2018hierarchical}; \cite{barrett2020evolution}). A hierarchical model is a collection of agents that work together to develop more sophisticated signalling systems. Applying the CF framework to such collections of agents to characterize precisely the ways in which they can combine to form more sophisticated agents could help us better understand existing hierarchical models, as well as suggest new directions of investigation. It might also help provide a formal justification for moving from models of language learning in which we consider a \textit{community} of learners (multiple agents) to one in we consider a single agent learning (who is composed of subagents working in concert). 

Additionally, we can use the subagency notion to model the same agent across time. In particular, we can model an agent before and after it makes a \textit{commitment} to do something. It is interesting that the same subagent notion that we can use to understand the relationships between different agents (as in the cases above) easily generalizes to account for this case (treating the agent with fewer options as a subagent of the earlier agent that hasn't yet made a commitment).\footnote{See \cref{comm} for the technical details of committing.}
\\
\\
We expect this framework to interface with many other topics in formal philosophy. Indeed, it is often true that new areas of investigation open up once we have a new vocabulary and formalism to express ideas. Thus, the primary goal of this paper is to introduce core parts of the technical apparatus to a philosophical audience, in the hopes that this will equip philosophers with tools they can use to tackle old and new problems. Having provided a sketch of how this framework might help with philosophical questions, we now turn to this goal, confident that readers will see the utility of the framework for themselves. In \cref{CF} and \cref{bi} we introduce essential elements of the framework. In \cref{coarse} we show how the framework provides a precise way to think about differing levels of description. In \cref{sub} we introduce the formal notion of subagency, and show how the framework provides a calculus of subagency.

\section{Technical Framework}
\label{tech}

\subsection{Cartesian Frames}
\label{CF}

To motivate the following technical definitions, consider one common way of representing decision problems.
Decision problems are often represented by matrices, where the rows represent a choice by the agent and the columns represent choices by the other agent(s) in a game, or the state of the environment.
For example, in a Newcomb problem, the two rows represent the two options available to the agent: either one-box or two-box.
The two columns represent the possible states of the environment---in this case, whether the opaque box contains money or does not.
Finally, the entries in the matrix represent the outcome of this interaction between the agent's choice and the state of the environment.
So we can think of a decision problem matrix as a combination of three things: (i) a set of options for the agent, (ii) a set of options for the environment,\footnote{Which can include other agents and their acts.} and (iii) the outcome produced by the combination of an option for the agent and an option for the environment.
We can generalize this matrix representation:

\bigskip
\noindent \textbf{Definition 1}.
Let $W$ be a set.
A Cartesian frame $C$ over $W$ is a triple $C = (A, E, \cdot)$, where $A, E$ are sets and $\cdot : A \times E \to W$ is a function.

\bigskip
\noindent The intended interpretation is: $W$ is a set, often possible worlds but occasionally utilities, which we can think of as possible outcomes. $A$ is the set of acts available to the agent, $E$ is the set of states of the environment, and given $a \in A, e \in E$, $a \cdot e$ is the possible world resulting from taking act $a$ in environment $e$.  
We can represent this idea more concretely with the following matrix:

\bigskip
\[
\begin{blockarray}{cccc}
& e_1 & e_2 & e_3\\
\begin{block}{c(ccc)}
a_1 & w_1 & w_2 & w_3\\
a_2 & w_4 & w_5 & w_6\\
a_3 & w_7 & w_8 & w_9\\
\end{block}
\end{blockarray}
\]

\bigskip
\noindent Here $A = \{ a_1, a_2, a_3\}$ represents an agent with three possible acts.
Similarly $E = \{e_1, e_2, e_3\}$ represents an environment with three possible states.
Finally, we represent the fact that the outcome, a possible world in $W$, is the result of the joint interaction of the state of the agent and the state of the environment by writing, e.g., $a_1 \cdot e_2 = w_2$.

Here we introduce a bit of notation.
If $C = (A, E, \cdot)$ is a Cartesian frame, we write $Agent(C)$ to denote $A$, the set of options for the agent;
we write $Env(C)$ to denote $E$, the set of options for the environment;
and we write $Image(C)$ to denote the set $\{w \in W : \exists a \in A, \exists e \in E (a \cdot e = w) \}$.
In other words, the image of a Cartesian frame is the subset of $W$ which appears in the matrix representation of that frame;
equivalently, the image is the subset of possible worlds that can be produced by some choice of agent and some choice of environment in that frame.

So a Cartesian frame can be considered as a generalized decision matrix.\footnote{Note that there are no cardinality restrictions in \textbf{Definition 1}; thus the framework allows infinite decision matrices as well.}
But we go a step further.
One of the novel strengths of Cartesian frames is the fact that one can define \emph{morphisms} between two frames, allowing a sort of transformation from one to the other:

\bigskip
\noindent \textbf{Definition 2}.
A \emph{morphism} from $C = (A, E, \cdot)$ to $D = (B, F, \star)$ is a pair of functions $(g: A \to B, h: F \to E)$ satisfying $a \cdot h(f) = g(a) \star f$ for all $a \in A, f \in F$.
We call this equation the \emph{adjointness condition} for morphisms.
For two morphisms $(g_0, h_0), (g_1, h_1)$, their composition is given by $(g_1, h_1) \circ (g_0, h_0) = (g_1 \circ g_0, h_0 \circ h_1)$.

\bigskip

Consider as an example a driver who is deciding whether to drive by the seaside ($a_{s}$) or whether to drive on the highway ($a_{h}$). There are three possible weather conditions: rainy ($e_{r}$), cloudy ($e_{c}$), and sunny ($e_{s}$). If there is sun and she takes the seaside route then she will get to her destination on time and she will enjoy a good view. If it is cloudy she will still get there on time, but the view will be obscured. If it is rainy she will be late and the view will be obscured. The highway is always fast, and so she will get to her destination on time. We represent her decision situation with the following Cartesian frame, where the worlds are numbers representing how much she likes each outcome.

\bigskip
\[
\begin{blockarray}{cccc}
& e_r & e_c & e_s\\
\begin{block}{c(ccc)}
a_s & 1 & 5 & 7\\
a_h & 5 & 5 & 5\\
\end{block}
\end{blockarray}
\]

Now imagine that the driver discovers a new road that goes through the countryside and has a nice view. Imagine that this road is just as fast as the highway, and is not affected by weather. We can represent this different set of options with the frame

\bigskip
\[
\begin{blockarray}{cccc}
& f_r & f_c & f_s\\
\begin{block}{c(ccc)}
b_s & 1 & 5 & 7\\
b_h & 5 & 5 & 5\\
b_c & 6 & 6 & 6\\
\end{block}
\end{blockarray}
\]

\noindent There is a morphism from the previous frame to this frame. (Consider the morphism $(g, h)$ where both $g$ and $h$ map elements to elements with corresponding subscripts.)  Morphisms are useful because they allow us to relate one Carteisan frame to another, and in this way define other useful concepts (many of which we will present in this paper). Note that in a morphism from $C$ to $D$, one function maps ``forward'' from the agent of $C$ to the agent of $D$, whereas the other function maps ``backward'' from the environment of $D$ to the environment of $C$.

Mathematically, Cartesian Frames are exactly Chu spaces.
Chu spaces provide a natural generalization of matrices, algebras, posets, and related algebraic structures.
Most often Chu spaces are studied from a categorical perspective.
Categorically, $\text{Chu}(W)$ is a category having as objects all Cartesian frames defined over $W$ and having as arrows all morphisms between those frames.
In the main body of this paper we will focus on the matrix representation of Cartesian frames, deferring categorical details to an appendix.\footnote{The mathematically inclined reader can find more information about Chu spaces at http://chu.stanford.edu/guide.html.}

\subsection{Biextensional Equivalence}
\label{bi}

One can also define multiple notions of \emph{equivalence} of Cartesian frames. These equivalence notions will be useful for us because they will allow us to abstract away inessential details of the Cartesian frames, helping us to focus on the decision-relevant structure. 
As in other parts of mathematics, we can define a notion of isomorphism:

\bigskip
\noindent \textbf{Definition 3}.
A morphism $(g, h) : C \to D$ is an \emph{isomorphism} if both $g$ and $h$ are bijections.
If there is an isomorphism between $C$ and $D$, we say $C \cong D$.

\bigskip
\noindent Isomorphism is a useful notion, but in fact we will need something more specific to Cartesian frames.
Intuitively, two Cartesian frames might differ and yet represent the ``same situation'' by containing duplicate rows and columns---that is, two or more rows or columns with exactly the same entries.
In a sense, if two rows (acts) result in exactly the same outcomes, then as far as the agent is concerned, they are the same act.
So if we have two CFs---one which has no duplicates, and another which is identical to the first except that it contains duplicates---then these two frames intuitively do not differ in the choices available to the agent.
Thus we want a way to identify two such Cartesian frames mathematically, despite the fact that they are not isomorphic.
We call this identification \emph{biextensional equivalence}.
We will define this relation in steps.
First, a Cartesian frame is \emph{biextensional} just in case it contains no ``duplicate'' rows and columns:

\bigskip
\noindent \textbf{Definition 4}.
A Cartesian frame $C= (A, E, \cdot)$ is called \emph{biextensional} if whenever $a_0, a_1 \in A$ are such that $a_0 \cdot e = a_1 \cdot e$ for all $e \in E$, then $a_0 = a_1$, and whenever $e_0, e_1 \in E$ are such that $a \cdot e_0 = a \cdot e_1$ for all $a \in A$, then $e_0 = e_1$.

\bigskip
\noindent Of course, many frames are not biextensional.
But from a given Cartesian frame $C$ we can generate a corresponding biextensional frame $\hat{C}$, called the \emph{biextensional collapse} of $C$.
We do this by ``deleting'' duplicate rows and columns:

\bigskip
\noindent \textbf{Definition 5}.
Given a Cartesian frame $C = (A, E, \cdot)$, for $a_0, a_1 \in A$, we say $a_0 \sim a_1$ if $a_0 \cdot e = a_1 \cdot e$ for all $e \in E$. For $e_0, e_1 \in E$, we say $e_0 \sim e_1$ if $a \cdot e_0 = a \cdot e_1$ for all $a \in A$.

Then for $a \in A$, let $\hat{a}$ denote the equivalence class of $a$ up to $\sim$. 
Let $\hat{A}$ denote the set of equivalence classes of $\sim$ in $A$.
Similarly, let $\hat{e}$ denote the equivalence class of $e$ up to $\sim$, and let $\hat{E}$ denote the set of equivalence classes of $\sim$ in $E$.

\bigskip
\noindent \textbf{Definition 6}.
Given a Cartesian frame $C = (A, E, \cdot)$, the \emph{biextensional collapse} of $C$ is given by $\hat{C} = (\hat{A}, \hat{E}, \hat{\cdot})$, where $\hat{a} \hat{\cdot} \hat{e} = a \cdot e$.

\bigskip
\noindent Putting it all together, we get:

\bigskip
\noindent \textbf{Definition 7}.
Given two Cartesian frames $C$ and $D$, we say that $C$ and $D$ are \emph{biextensionally equivalent}, written $C \simeq D$, just in case $\hat{C} \cong \hat{D}$.

\bigskip
\noindent \textbf{Claim 8} \label{claim:isothenbiex}.
For two Cartesian frames $C$ and $D$, if $C = D$ then $C \cong D$, and if $C \cong D$ then $C \simeq D$.

\begin{proof}
If $C = D$, then $(\text{id}_A, \text{id}_E)$ is an isomorphism between $C$ and $D$.

For the proof of the second claim, see Appendix A.
\end{proof}

\noindent We now have a precise way to express the idea that Cartesian frames with duplicate rows/columns express ``the same'' situation.
Two frames might be nonisomorphic, yet if duplicates are deleted, the resulting biextensional collapses might be isomorphic. For example, the following two frames are not isomorphic, but they are biextensionally equivalent (consider deleting $b_{3}$).

\bigskip
\[
\begin{blockarray}{cccc}
& e_1 & e_2 & e_3\\
\begin{block}{c(ccc)}
a_1 & w_1 & w_2 & w_3\\
a_2 & w_4 & w_5 & w_6\\
\end{block}
\end{blockarray}
\hspace{2cm}
\begin{blockarray}{cccc}
& f_1 & f_2 & f_3\\
\begin{block}{c(ccc)}
b_1 & w_1 & w_2 & w_3\\
b_2 & w_4 & w_5 & w_6\\
b_3 & w_4 & w_5 & w_6\\
\end{block}
\end{blockarray}
\]

\noindent In general we will be working up to biextensional equivalence.\footnote{There is another handy equivalence relation on Cartesian frames called ``homotopy equivalence''. It is in fact provably equivalent to biextensional equivalence. We will not make use of it until Appendix A.}

We also introduce the \emph{dual} of a Cartesian frame, which will be important in what follows:

\bigskip
\noindent \textbf{Definition 9}.
Let $-^* : \text{Chu}(W) \to \text{Chu}(W)^{\text{op}}$ be the functor given by $(A, E, \cdot)^* = (E, A, \star)$, where $e \star a = a \cdot e$.\footnotemark
\footnotetext{Note that here we mean `functor' in the categorical sense, i.e., a map from one category to another which preserves domains and codomains of arrows, preserves composition of arrows, and maps identities to identities.}

\bigskip
\noindent Note that, if we consider a given Cartesian frame as a matrix, $-^*$ is just transposition.

\subsection{Coarse and Refined Worlds}
\label{coarse}

Cartesian frames are defined over a set of possible worlds, say, $V$.
We can think of these possible worlds as carving up the world at a certain level of description.
Perhaps $V$ has separate possible worlds for the scenario where Amy has a sandwich for lunch versus the scenario where Amy has spaghetti.
But one could easily imagine that the worlds in $V$ do not distinguish, say, the scenario where Amy has a sandwich \emph{and} wears red socks versus the scenario where Amy has a sandwich and wears yellow socks.
Perhaps there's only one world in $V$ in which Amy has a sandwich at all, and thus $V$ cannot distinguish the above scenarios.
But one could instead consider a set of possible worlds $W$ which \emph{does} include possible worlds that distinguish these scenarios.
Thus, at least with respect to Amy's socks, $W$ embodies a ``finer-grained'' level of description of the world.
Cartesian frames offer a natural mechanism for moving between such levels of description.

\bigskip
\noindent \textbf{Definition 10}.
Given two sets $W$ and $V$, and a function $p: W \to V$, let $p^{\circ} : Chu(W) \to Chu(V)$ denote the functor that sends the object $(A, E, \cdot) \in Chu(W)$ to $(B, F, \star) \in Chu(V)$, where $b \star f = p(a \cdot e)$, and sends the morphism $(g, h)$ to the morphism with the same underlying functions, $(g, h)$.\footnote{Note that $p^{\circ}$ preserves biextensional equivalence, i.e., if $C \simeq D$ then $p^{\circ}(C) \simeq p^{\circ}(D)$.}

\bigskip
\noindent These maps are \emph{functors} between categories, hence the name.

If $p$ is surjective, then we say that ``$W$ is a refined version of $V$'', or ``$V$ is a coarse version of $W$''.
To see why, imagine representing the set $W$ of possible worlds as a two-dimensional plane.
Each point is a specific possible world.
But we can draw boundaries around groups of worlds, thereby assigning worlds to ``neighborhoods''.
We then treat the neighborhoods as possible worlds themselves, and group them into a new set, $V$.
Thus every world in $W$ gets mapped to a world in $V$, i.e., the function is surjective.
Some propositions are true (false) at each world in the neighborhood;
these propositions are true (false) at the new neighborhood-world.
Other propositions are true at some worlds in the neighborhood and false at others;
these propositions are thus not settled by the neighborhood-world.
In other words, the neighborhood-world is ``less detailed'' than the worlds it contains, or, equivalently, is a ``coarse version'' of those worlds.\footnotemark
\footnotetext{There's a clear similarity here between coarsening/refining functions between sets of worlds and the coarsening/refining relation found in possibility semantics (originally proposed in \cite{humberstone1981possibilities}, with substantial development in \cite{holliday2018possibilityframes}.)}

For example, let $W=\{a, o, c\}$ and $V=\{2, 3\}$, and let $p(a)=p(o)=2, p(c)=3$. We can think of $a, o$, and $c$ as apple, orange, and chocolate respectively, with the elements of $V$ representing subjective levels of enjoyment of the agent. We can then think of $p$ as encoding how much an agent enjoys a particular snack. Now consider the following Cartesian frame, $D$, over $W$

\bigskip
\[
\begin{blockarray}{ccccc}
& O & A & OC & AC\\
\begin{block}{c(cccc)}
C & o & a & c & c\\
\neg C & o & a & o & a\\
\end{block}
\end{blockarray}
\]

\noindent where $O$ and $A$ are worlds in which there are only oranges and apples available to the agent respectively, $OC$ is a world in which both oranges and chocolate are available, and $AC$ is a world in which both apples and chocolate are available. The agent has two options, asking for chocolate ($C$), and not asking for chocolate ($\neg C$). Clearly $D$ is biextensional. Now consider $p^{\circ}(D)$: 

\bigskip
\[
\begin{blockarray}{ccccc}
& Or & Ap & OCh & ACh\\
\begin{block}{c(cccc)}
Ch & 2 & 2 & 3 & 3\\
\neg Ch & 2 & 2 & 2 & 2\\
\end{block}
\end{blockarray} 
\hspace{0.5cm}
\simeq 
\hspace{0.5cm}
\begin{blockarray}{ccc}
& Fr & FrCh \\
\begin{block}{c(cc)}
Ch & 2 & 3\\
\neg Ch & 2 & 2\\
\end{block}
\end{blockarray}\ 
\]

Notice how $p^{\circ}(D)$ loses the information about which fruit the agent gets, and instead the worlds only describe how much the agent values their snack. Furthermore, notice how once we coarsen the situation like this we can drop decision-irrelevant environments by taking the biextensional collapse to yield the Cartesian frame on the right. Thus we see how the CF framework gives us an easy way to non-trivially move between decision problems at different levels of description.

Earlier we presented a surjective function $p: W \to V$, and called the codomain a ``coarsening''. Given a set $W$ and a coarsening, $V$, there is one unique function from $W$ to $V$, namely the function that assigns each world to its ``neighborhood''. But the reverse direction is not unique; given a neighborhood, there are many worlds to which a function could map.

We solve this problem by introducing functors defined from Cartesian frames, rather than defined from functions on sets of possible worlds.
Intuitively, when we define a functor from a Cartesian frame $C$, the agent of $C$ chooses a neighborhood in $V$, and the environment of $C$ chooses a function that selects an element in that neighborhood, i.e., the environment chooses a function $q : V \to W$.

\bigskip
\noindent \textbf{Definition 11}.
Let $C = (V, E, \cdot)$ be a Cartesian frame over $W$, whose agent is $V$.
Then $C^{\circ} : \text{Chu}(V) \to \text{Chu}(W)$ is the functor that sends $(B, F, \star)$ to $(B, F \times E, \diamond)$, where $b \diamond (f, e) = (b \star f) \cdot e$, and sends the morphism $(g, h)$ to $(g, h')$ where $h'(f, e) = (h(f), e)$.

\bigskip
\noindent Note that this definition looks a bit like ``currying'', the idea that there is a natural correspondence between functions of $n+1$ arguments and functions from $n$ arguments to a function space. For example, if $f: A \times B \to C$, there is a corresponding ``curried'' function $\hat{f} : A \to C^B$.
In our example we see a natural correspondence between ``functions'' $k: B \times F \to W^E$ (i.e., $(B, F, \star)$ is defined over $V$, the agent of $C$, which one could consider as a set of functions from $E$ to $W$) and $\hat{k} : B \times F \times E \to W$.

\subsection{Subagents}
\label{sub}

We arrive at the (in our opinion) most intriguing feature of Cartesian frames: subagency.\footnote{Recall our informal discussion of subagency in \cref{back}.} In order to give the categorical definition of a subagent, we introduce one additional definition. 

\bigskip
\noindent \textbf{Definition 12}.
Given $S \subseteq W$, $\bot_S$ is the Cartesian frame $\bot_S = (S, \{e\}, \star)$, where $\star$ is given by $s \star e = s$ for all $s \in S$.
$\bot$ is the Cartesian frame $\bot_W$.

\bigskip
\noindent \textbf{Definition 13} (Categorical definition of subagent).
Let $C$ and $D$ be Cartesian frames over $W$.
We say that $C$ is a subagent of $D$, written $C \triangleleft D$, if for every morphism $\phi : C \to \bot$ there exists a pair of morphisms $\phi_0 : C \to D$, $\phi_1 : D \to \bot$ such that $\phi = \phi_1 \circ \phi_0$ (i.e., $\phi$ factors through $D$).
In other words, the following diagram commutes:

\begin{center}
 \begin{tikzcd}[row sep=large, column sep=large]
C \arrow[r, "\phi"] \arrow[d, "\phi_0"] & \bot \\ 
D \arrow[ru, "\phi_1", swap]
\end{tikzcd}
\end{center}

\noindent To get a feel for the significance of this definition, consider the fact that the morphisms $\phi_e : C \to \bot$ correspond exactly to elements of $E = Env(C)$.
To see this, let $\bot = (W, \{j\}, \star)$.
Let the morphism $(g, h) : C \to \bot$ be given by $g(a) = a \cdot e$, and $h(j) = e$.
This morphism satisfies the adjointness condition:
\begin{displaymath}
\begin{split}
    a \cdot h(j) & = a \cdot e \\
                 & = w \\
                 & = w \star j \\
                 & = g(a) \star j
\end{split}
\end{displaymath}
So the categorical definition states that every element of $Env(C)$ is the image of some element of $Env(D)$ under some morphism $\phi_0 : C \to D$.

It will sometimes be useful for us to work with a different definition of subagency.

\bigskip 

\noindent \textbf{Definition 14}.
(Currying definition of subagent).
Let $C$ and $D$ be Cartesian frames over $W$.
We say that $C \triangleleft D$ if there exists a Cartesian frame $Z$ over Agent($D$) such that $C \simeq D^{\circ}(Z)$.
\bigskip

\noindent \textbf{Claim 15}.
The categorical definition and the currying definition of subagency are equivalent. 

\begin{proof}
See Appendix B.
\end{proof}

\bigskip
\noindent We can think of $D^{\circ}$ and $Z$ as ``factoring'' $C$'s agent into a further interaction of an agent and an environment.
Consider:

\[
\hspace{-0.5cm} 
Z
\hspace{3cm}
D
\hspace{3.5cm}
D^{\circ}(Z) \simeq C
\]
\[
\begin{blockarray}{ccc}
& f_1 & f_2 \\
\begin{block}{c(cc)}
b_1 & v_1 & v_1 \\
b_2 & v_2 & v_3 \\
\end{block}
\end{blockarray}
\hspace{1cm}
\begin{blockarray}{ccc}
& e_1 & e_2 \\
\begin{block}{c(cc)}
v_1 & w_0 & w_1 \\
v_2 & w_2 & w_3 \\
v_3 & w_4 & w_5 \\
\end{block}
\end{blockarray}
\hspace{1cm}
\begin{blockarray}{ccccc}
& f_{1}e_1 & f_{1}e_2 & f_{2}e_1 & f_{2}e_2 \\
\begin{block}{c(cccc)}
b_1 & w_0 & w_1 & w_0 & w_1 \\
b_2 & w_2 & w_3 & w_4 & w_5 \\
\end{block}
\end{blockarray}
\]

\noindent In this case, $D$ is a frame whose agent is $V$ and whose environment is $E$.
$Z$ is a frame which represents a ``factoring'' of the agent of $D$;
we see that it is the interaction of an agent, $B$, and an environment, $F$, and their joint interaction results in a choice $v_n \in V$.
We can liken $Z$ to the perspective of a single player, $B$, in a team, $V$.
The behavior of $V$ is jointly determined by $B$'s behavior and the behavior of the rest of the team, $F$.

So $Z$ represents the interaction of a single player with the rest of the team.
But we can ``zoom out'' and consider the team's interaction with its environment, $E$, \emph{from the perspective of $B$}.
That is, we can consider a frame whose agent is $B$ and whose environment is made up of both (i) choices made by the rest of the team and (ii) states of the team's environment.
This is the action of $D^{\circ}$ on $Z$.
Thus $D^{\circ}(Z)$ represents the perspective of the single player interacting with both their team and their team's environment. This frame is biextensionally equivalent to $C$, so $C$ can likewise be thought of as the perspective of a single player in the original team.
It is in this sense that we can claim that $C$ is a subagent of $D$.

\bigskip
\noindent \textbf{Claim 16}.
$\triangleleft$ is transitive.

\begin{proof}
If $C_0 \triangleleft C_1$ and $C_1 \triangleleft C_2$, then, by the categorical definition, every morphism $\phi_0 : C_0 \to \bot$ factors through $C_1$ and every morphism $\phi_1 : C_1 \to \bot$ factors through $C_2$.
Thus the following diagram commutes:

\begin{center}
    \begin{tikzcd}[row sep=large, column sep=large]
    C_0 \arrow[r, "\phi_0"] \arrow[d] & \bot \\
    C_1 \arrow[ru, "\phi_1"] \arrow[r, swap] & C_2 \arrow[u]
    \end{tikzcd}
\end{center}
So $\phi_0$ factors through $C_2$.
\end{proof}

\bigskip
\noindent \textbf{Claim 17}.
$\triangleleft$ is reflexive. Further, if $C \simeq D$, then $C \triangleleft D$.

\begin{proof}
Let $C = (A, E, \cdot)$ and $D = (B, F, \star)$ be Cartesian frames over $W$, with $C \simeq D$.
Consider the Cartesian frame $Z$ over $B$ given by $Z = (B, \{x\}, \diamond)$, where $b \diamond x = b$.
Then $D^{\circ}(Z) = (B, \{x\} \times F, \bullet)$.
Observe that $D \cong D^{\circ}(Z)$, so $C \triangleleft D$, according to the currying definition.

Clearly $C \simeq C$, so $C \triangleleft C$.
\end{proof}

\noindent Our notion of subagency can be further refined into two distinct notions: additive and multiplicative subagents.
Intuitively, an additive subagent is an agent who faces the same environment but has fewer options.
We call such a subagent ``additive'' because one must simply add more options to the subagent to return to the original superagent.
A multiplicative subagent is better embodied by the relationship between a team and an individual member.
We call such a subagent ``multiplicative'' because members act \emph{jointly} to produce the action of the team;
in this sense, an action of a team can be ``factored'' into some tuple of actions by the individual members, or equivalently, the team is the Cartesian product of the members. We present two pairs of (equivalent) definitions of additive and multiplicative subagent.

The first definition of additive subagency is what we call the \emph{committing} definition.
Intuitively, we are viewing $C$ as the result of $D$ making a commitment to only select from some subset of their available acts. 

\bigskip
\noindent \textbf{Definition 18} (Committing definition).
Given Cartesian frames $C$ and $D$ over $W$, we say $C \triangleleft_{+} D$, read ``$C$ is an additive subagent of $D$'', if there exist three sets $X$, $Y$, and $Z$, with $X \subseteq Y$, and a function $f : Y \times Z \to W$ such that $C \simeq (X, Z, \diamond)$ and $D \simeq (Y, Z, \bullet)$, where $\diamond$ and $\bullet$ are given by $x \diamond z = f(x, z)$ and $y \bullet z = f(y, z)$.

\bigskip
\noindent As an example, consider our earlier example of the driver with the following two Cartesian frames:
\[
\begin{blockarray}{cccc}
& e_r & e_c & e_s\\
\begin{block}{c(ccc)}
a_s & 1 & 5 & 7\\
a_h & 5 & 5 & 5\\
\end{block}
\end{blockarray}
\hspace{2cm}
\begin{blockarray}{cccc}
& f_r & f_c & f_s\\
\begin{block}{c(ccc)}
b_s & 1 & 5 & 7\\
b_h & 5 & 5 & 5\\
b_c & 6 & 6 & 6\\
\end{block}
\end{blockarray}
\]

\noindent The agent of the frame on the left is an additive subagent of the agent of the frame on the right. We can think of the left frame as arising from the agent on the right committing to not taking the country road. 

\bigskip
\noindent We next give our first definition of multiplicative subagency.
We refer to this definition as the \emph{externalizing} definition of multiplicative subagency.
Intuitively, we are viewing $C$ as the result of $D$ sending some of its choices into the environment.
For example, if the agent of $D$ is a team, then we can send some members into the environment to obtain the subagent frame which represents the remaining members.

\bigskip
\noindent \textbf{Definition 19} (Externalizing definition).
Given Cartesian frames $C$ and $D$ over $W$, we say $C \triangleleft_{\times} D$, read ``$C$ is a multiplicative subagent of $D$'', if there exist three sets $X, Y,$ and $Z$, and a function $f : X \times Y \times Z \to W$ such that $C \simeq (X, Y \times Z, \diamond)$ and $D \simeq (X \times Y, Z, \bullet)$, where $\diamond$ and $\bullet$ are given by $x \diamond (y, z) = f(x, y, z)$ and $(x, y) \bullet z = f(x, y, z)$.

\bigskip
\noindent Consider our earlier example from the discussion of the currying definition of subagency.
$D^{\circ}(Z)$ is in fact a multiplicative subagent of $D$.
We see that $D$'s actions can be factored into members of $B$ and $F$, i.e., the joint interaction of $D^{\circ}(Z)$'s agent and some other subagent of $D$;
we relabel $v_i \in V$ to reflect the fact that they can be factored in this manner.

\[
Z
\hspace{3cm}
D
\hspace{3cm}
D^{\circ}(Z) \simeq C
\]
\[
\begin{blockarray}{ccc}
& f_1 & f_2 \\
\begin{block}{c(cc)}
b_1 & v_1 & v_1 \\
b_2 & v_2 & v_3 \\
\end{block}
\end{blockarray}
\hspace{.5cm}
\begin{blockarray}{ccc}
& e_1 & e_2 \\
\begin{block}{c(cc)}
b_{1}f_1=v_1 & w_0 & w_1 \\
b_{2}f_1=v_2 & w_2 & w_3 \\
b_{2}f_2=v_3 & w_4 & w_5 \\
\end{block}
\end{blockarray}
\hspace{.5cm}
\begin{blockarray}{ccccc}
& f_{1}e_1 & f_{1}e_2 & f_{2}e_1 & f_{2}e_2 \\
\begin{block}{c(cccc)}
b_1 & w_0 & w_1 & w_0 & w_1 \\
b_2 & w_2 & w_3 & w_4 & w_5 \\
\end{block}
\end{blockarray}
\]

\bigskip
\noindent We now present another pair of definitions of additive and multiplicative subagents, equivalent to the first pair.
While the first pair of definitions emphasizes the relationships between the sets of agent choices and the set of environment states, this new pair will rely on the notion of ``currying'' we developed via functors from Cartesian frames.
We therefore refer to the following definitions as the \emph{currying} definitions of subagency.

\bigskip
\noindent \textbf{Definition 20} (Currying definition of additive subagent).
We say that $C \triangleleft_{+} D$ if there exists a Cartesian frame $M$ over $Agent(D)$ with $|Env(M)| = 1$, such that $C \simeq D^{\circ}(M)$.

\bigskip
\noindent \textbf{Definition 21} (Currying definition of multiplicative subagent).
We say that $C \triangleleft_{\times} D$ if there exists a Cartesian frame $M$ over $Agent(D)$ with $Image(M) = Agent(D)$ such that $C \simeq D^{\circ}(M)$. 

\bigskip
\noindent \textbf{Claim 22}. \label{claim:subcurrycommequiv}
The currying definitions of additive subagent is equivalent to the committing definition, and the currying definition of multiplicative subagent is equivalent to the externalizing definition.

\begin{proof}
See Appendix A.
\end{proof}

\noindent We mention here a few nice properties of the additive and multiplicative subagent relations:

\bigskip
\noindent \textbf{Claim 23}.
\begin{enumerate}
    \item If $C \triangleleft_{+} D$, then $C \triangleleft D$. If $C \triangleleft_{\times} D$, then $C \triangleleft D$.
    \item If $C \triangleleft_{+} D$, $C \simeq C'$, and $D \simeq D'$, then $C' \triangleleft_{+} D'$. Similarly, if $C \triangleleft_{\times} D$, $C \simeq C'$, and $D \simeq D'$, then $C' \triangleleft_{\times} D'$.
    \item Both $\triangleleft_{+}$ and $\triangleleft_{\times}$ are reflexive and transitive.
\end{enumerate}

\noindent We are now in a position to prove our first theorem:

\bigskip
\noindent \textbf{Theorem 24} (Decomposition Theorem).
$C_0 \triangleleft C_1$ iff there exists a $C_2$ such that $C_0 \triangleleft_{\times} C_2 \triangleleft_{+} C_1$.

\begin{proof}
We will use the currying definitions of subagent and multiplicative subagent, and the committing definition of additive subagent.
Let $C_0 = (A_0, E_0, \cdot_0)$ and $C_1 = (A_1, E_1, \cdot_1)$.
If $C_0 \triangleleft C_1$, then there exists some Cartesian frame $D$ over $A_1$ such that $C_0 = C^{\circ}_1 (D)$. 

Let $C_2 = (Image(D), E_1, \cdot_2)$, where $\cdot_2$ is given by $a \cdot_2 e = a \cdot_1 e$.
Since $Image(D) \subseteq A_1$, $C_2$ is created by deleting some possible agents from $C_1$, so by the committing definition of additive subagent, $C_2 \triangleleft_+ C_1$.

By the currying definition of multiplicative subagent, we want to show that there is a Cartesian frame $M$  over $Agent(C_2)$ with $Image(M) = Agent(C_2)$ such that $C_0 \simeq C^{\circ}_2 (M)$.
But $C^{\circ}_1 (D) = C^{\circ}_2 (D)$, since $C_1$ and $C_2$ have the same environment.
Thus $C_0 = C^{\circ}_2 (D)$, so $C_0 \simeq C^{\circ}_2 (D)$, and $Image(D) = Agent(C_2)$, so $C_0 \triangleleft_{\times} C_2$.
\end{proof}

This theorem shows that if one agent is a subagent of another, then there is always a way to break down that subagency relationship into (possibly) multiple subagency relationships that are all multiplicative and additive. That is, we can always gain clarity about a certain subagency relationship by expressing it in terms of the more familiar multiplicative and additive relations. 

\subsubsection{Committing, Assuming, Externalizing, Internalizing}
\label{comm}

Our first pair of definitions of additive and multiplicative subagents were called ``committing'' and ``externalizing'', to reflect the idea that one can generate subagents from a CF in particular ways.
We now define these operations and their inverses more rigorously, providing us with a set of tools for decomposing agents into subagents and composing agents into superagents.
First, we introduce the notion of a sub-\emph{environment}, dual to subagent:

\bigskip
\noindent \textbf{Definition 25} (Sub-environment).
We say $C$ is a sub-environment of $D$, written $C \triangleleft^* D$, if $D^* \triangleleft C^*$.

\bigskip
\noindent The notion of a sub-environment can naturally be refined into additive and multiplicative notions:

\bigskip
\noindent \textbf{Definition 26}.
We say $C$ is an additive sub-environment of $D$, written $C \triangleleft^*_{+} D$, if $D^* \triangleleft_+ C^*$.
Similarly, we say $C$ is a multiplicative sub-environment of $D$, written $C \triangleleft^*_{\times} D$, if $D^* \triangleleft_{\times} C^*$.

\bigskip
\noindent In fact, subagent and sub-environment coincide when the notion is multiplicative:

\bigskip
\noindent \textbf{Claim 27}.
$C \triangleleft_{\times} D$ iff $C \triangleleft^*_{\times} D$.

\begin{proof}
Assuming the externalizing definition of $\triangleleft_{\times}$, we have $D^* \simeq (Z, Y \times X, \cdot)$ and $C^* \simeq (Z \times Y, X, \bullet)$.
This is exactly the externalizing definition for $D^* \triangleleft_{\times} C^*$, so $C \triangleleft^*_{\times} D$.

Conversely, if $C \triangleleft^*_{\times} D$, then $D^* \triangleleft_{\times} C^*$, so $C \cong (C^*)^* \triangleleft_{\times} (D^*)^* \cong D$.
\end{proof}

\bigskip
\noindent \textbf{Definition 28} (Committing).
Given a Cartesian frame $C = (A, E, \cdot)$ and a subset $B \subseteq A$, let $Commit^B (C)$ denote the Cartesian frame $(B, E, \star)$ with $\star$ given by $b \star e = b \cdot e$.
Let $Commit^{\setminus B} (C)$ denote the Cartesian frame $(A \setminus B, E, \diamond)$ with $\diamond$ given by $a \diamond e = a \cdot e$.

\bigskip
\noindent $Commit^B$ represents the Cartesian frame which results from the agent of $C$ committing to choose an element of $B$ (\emph{mutatis mutandis} for $\setminus B$).
In other words, the agent reduces their available acts to some subset of their original available acts.

\bigskip
\noindent \textbf{Definition 29} (Assuming).
Given a subset $F \subseteq E$, let $Assume^F (C)$ denote the Cartesian frame $(A, F, \star)$, where $\star$ is given by $a \star f = a \cdot f$. Let $Assume^{\setminus F}$ denote the Cartesian frame $(A, E \setminus F, \diamond)$, where $a \diamond e = a \cdot e$.

\bigskip
\noindent $Assume^F$ represents the Cartesian frame which results from assuming that the environment is chosen from $F$ (\emph{mutatis mutandis} for $\setminus F$).

We see that committing and assuming do indeed generate additive subagents and additive super-environments, respectively:

\bigskip
\noindent \textbf{Claim 30}.
\begin{enumerate}[label=(\arabic*)]
    \item For all $B \subseteq A$, $Commit^B (C) \triangleleft_+ C$ and  $Commit^{\setminus B} (C) \triangleleft_+ C$.
    \item For all $F \subseteq E$, $C \triangleleft^{*}_+ Assume^F (C)$ and $C \triangleleft^{*}_+ Assume^{\setminus F} (C)$.
\end{enumerate}

\begin{proof}
(1) is trivial from the committing definition of additive subagent.
(2) is trivial from the committing definition together with the definition of sub-environment.
\end{proof}

\noindent We also have two operations for generating multiplicative subagents and super-environments.

\bigskip
\noindent \textbf{Definition 31}.
Given a partition $X$ of $Y$, let $Y/X$ denote the set of all functions $q$ from $X$ to $Y$ such that $q(x) \in x$ for all $x \in X$.

\bigskip
\noindent \textbf{Definition 32} (Externalizing).
Given a partition $B$ of $A$, let $External^B (C)$ denote the Cartesian frame $(A/B, B \times E, \star)$, where $\star$ is given by $q \star (b, e) = q(b) \cdot e$.
Let $External^{/B} (C)$ denote the Cartesian frame $(B, A/B \times E, \diamond)$ where $\diamond$ is given by $b \diamond (q, e) = q(b) \cdot e$.

\bigskip
\noindent Here the partition, $B$, and the set of functions $A/B$ represent two multiplicative subagents of $A$.
To see why, consider a two-person team, comprised of Players 1 and 2.
Player 1 has some available acts, which grant Player 1 some agency over which outcome the team effects.
But Player 1 does not \emph{completely} determine the team's behavior, since the team equally requires Player 2's input.
So in a sense, Player 1 can at most effect some set of possible team behaviors;
they can reduce the set of possible outcomes to the subset consistent with their selected action.
For this reason it makes sense to represent Player 1 as a \emph{partition};
they can only partially determine the team's act.

Then when Player 2 acts, the outcome space is reduced from the particular subset which Player 1 selects to a particular outcome;
in this way it makes sense to represent Player 2 as a \emph{function} from Player 1's partition to the set of possible worlds.

We now define the inverse operation, \emph{internalizing}:

\bigskip
\noindent \textbf{Definition 33} (Internalizing).
Given a partition $F$ of $E$, let $Internal^F (C)$ denote the Cartesian frame $(F \times A, E/F, \star)$ where $\star$ is given by $(f, a) \star q = a \cdot q(f)$
Let $Internal^{/F} (C)$ denote the Cartesian frame $(E/F \times A, F, \diamond)$ where $\diamond$ is given by $(q, a) \diamond f = a \cdot q(f)$. 

\bigskip
\noindent $Internal^F (C)$ represents the Cartesian frame which results from factoring the environment into equivalence classes, selecting an element of each equivalence class, and then sending the choice of a particular equivalence class into the agent.

Similarly, $Internal^{/F} (C)$ represents the Cartesian frame which results from the same factorization, but the choice of an element from each equivalence class is sent into the agent.

We see that externalizing and internalizing do indeed generate multiplicative subagents and multiplicative superagents, respectively:

\bigskip
\noindent \textbf{Claim 34}. \label{claim:externintern}
\begin{itemize}
    \item For all partitions $B$ of $A$, $External^B (C) \triangleleft_{\times} C$ and $External^{/B} (C) \triangleleft_{\times} C$.
    \item For all partitions $F$ of $E$, $C \triangleleft_{\times} Internal^F (C)$ and $C \triangleleft_{\times} Internal^{/F} (C)$.
\end{itemize}

\begin{proof}
See Appendix A.
\end{proof}

\noindent We note in passing that these four operations are well-defined up to biextensional equivalence.
Further, all four operations are idempotent:

\bigskip
\noindent \textbf{Claim 35}.
For any subset $B \subseteq A$, and for any partition $F$ of $E$,
\begin{itemize}
    \item $Commit^B$, $Commit^{\setminus B}$, $Assume^B$, and $Assume^{\setminus B}$ are idempotent.
    \item $External^F$, $External^{/F}$, $Internal^F$, and $Internal^{/F}$ are idempotent.
\end{itemize}

\noindent Idempotence makes sense for these operations. For example, suppose that you have some frame $C$ and externalize some choice to the environment, producing frame $C'$; if you try to externalize the \textit{same} choice, you should just get back $C'$, since in that frame the choice is already external. The same intuition holds for committing, assuming, and internalizing. 

\section{Conclusion}

We have presented the CF framework because we believe it can become a valuable tool for formal philosophy.
Of course, there is plenty of work left to be done. In addition to the avenues of philosophical investigation we discussed in \cref{back}, we mention a few others.

It would be useful to understand how the introduction of a probability measure over environments interacts with the rest of the framework. Probability would allow us to express claims like ``an agent C chooses an action that maximizes expected utility" within the framework.\footnote{Similarly, it would be good to integrate a theory of utility into the framework, though this seems more straightforward to do. The example in \cref{coarse} with the fruit already points in this direction.}  

It would also be valuable to get a better sense of the structure of the algebra of subagency. Is there a way to generate a set of subagents for a certain frame (or set of frames) that has desirable closure properties? If so, what kind of structure does this yield? 

On the more technical side, it would be interesting to develop the categorical framework into an independent foundation for Cartesian frames, and see to what extent this helps streamline proofs and helps to pose and solve problems in formal philosophy. Furthermore, there is a deep connection between linear logic and the theory of Cartesian frames, since one can define binary operations on frames that correspond in a precise way to those of linear logic. It would be interesting to see if this connection could prove fruitful for investigations of the CF framework.

\section*{Acknowledgements}

We would like to thank Alex Appel, Rob Bensinger, Tsvi Benson-Tilsen, Andrew Critch, Abram Demski, Sam Eisenstat, David Girardo, Evan Hubinger, Edward Kmett, Ramana Kumar, Alexander Gietelink Oldenziel, Steve Rayhawk, Nate Soares, Brian Skyrms, Simon Huttegger, Jeffrey Barrett, Cailin O'Connor, James Weatherall, Toby Meadows, Aydin Mohseni, and Claire Wang for helpful discussion and comments on earlier drafts. We also want to thank the Formal Social Epistemology research group at the University of California, Irvine for a useful discussion of our paper.

This research was supported by the Machine Intelligence Research Institute.

\begin{flushright}

  Scott Garrabrant\\
\emph{
  Machine Intelligence Research Institute\\
  Reno, NV, USA\\
  scott@intelligence.org\\
}
\vspace{0.5cm}
  Daniel A. Herrmann\\
\emph{
  Department of Logic and Philosophy of Science\\
  University of California, Irvine\\
  Irvine, CA, USA\\
  daherrma@uci.edu\\
}
\vspace{0.5cm}
  Josiah Lopez-Wild\\
\emph{
  Department of Logic and Philosophy of Science\\
  University of California, Irvine\\
  Irvine, CA, USA\\
  jlopezwi@uci.edu\\
}
\end{flushright}

\appendix 

\section{Proofs of Claims}

Earlier in the paper we mentioned homotopy equivalence.
It will be useful in what follows, so we now define it and prove its equivalence to biextensional equivalence.
Since homotopy and biextensional equivalences are equivalent notions, we use `$\simeq$' to denote both (it will be clear from context which definition is in use).

\bigskip
\noindent \textbf{Definition 36} (Homotopic).
Two morphisms $(g_0, h_0), (g_1, h_1) : C \to D$ are \emph{homotopic} if $(g_0, h_1)$ is also a morphism.

\bigskip
\noindent One can think of $(g_0, h_1)$ as a sort of ``continuous deformation'' from $(g_0, h_0)$ to $(g_1, h_1)$, hence the term `homotopic'.

\bigskip
\noindent \textbf{Definition 37} (Homotopy equivalence).
$C$ is \emph{homotopy equivalent} to $D$, written $C \simeq D$, if there exists a pair of morphisms $\phi : C \to D$, $\psi : D \to C$, such that $\psi \circ \phi$ is homotopic to the identity on $C$, and $\phi \circ \psi$ is homotopic to the identity on $D$.

One can view homotopy equivalence as a weakening of isomorphism.
If $\phi$ and $\psi$ composed to \emph{equal} the identity in both orders, then $C$ and $D$ would be isomorphic.
Instead we only require that their compositions are homotopic to the identity in both orders.
As it turns out, homotopy equivalence and isomorphism coincide for biextensional frames.

\bigskip
\noindent \textbf{Claim 38}.
If $C = (A, E, \cdot)$ and $D = (B, F, \star)$ are biextensional, then $C \simeq D$ iff $C \cong D$.

\begin{proof}
``$\Leftarrow$''
If $C \cong D$, then there is some morphism $(g, h) : C \to D$ such that both of $g, h$ are bijective.
Then both of $g, h$ have inverses $g^{-1}, h^{-1}$.
Then $(g, h)$ and $(g^{-1}, h^{-1})$ compose to the identity in both orders, and $C \simeq D$.

``$\Rightarrow$''
Assume $C \simeq D$.
Then there is some $\phi : C \to D$ and some $\psi : D \to C$ that compose to something homotopic to the identity in both orders.
We will show that $\phi$ (and thus $\psi$) is an isomorphism.
Let $\phi = (g_{\phi}, h_{\phi})$ and $\psi = (g_{\psi}, h_{\psi})$.
Then $(g_{\psi} \circ g_{\phi}, \text{id}_E)$ is a morphism.
That is, 
\begin{displaymath}
g_{\psi} (g_{\phi} (a)) \cdot e =  a \cdot e
\end{displaymath}
for all $a \in A, e \in E$.
If for some $a$, $g_{\psi} (g_{\phi} (a)) \neq a$, then $A$ has two identical rows, i.e., $C$ is not biextensional.
So $g_{\psi} \circ g_{\phi} = \text{id}_A$, and $g_{\psi}, g_{\phi}$ are bijective.

We now show that $h_{\phi}$ is bijective.
$h_{\phi}$ is injective, since if $h_{\phi}(f_1) = h_{\phi}(f_2)$ then for all $a \in A$,
\begin{displaymath}
\begin{split}
    g_{\phi}(a) \star f_1 & = a \cdot h_{\phi}(f_1) \\
                          & = a \cdot h_{\phi}(f_2) \\
                          & = g_{\phi}(a) \star f_2
\end{split}
\end{displaymath}
So $f_1$ and $f_2$ are identical columns;
but $D$ is biextensional, so $f_1 = f_2$.

$h_{\phi}$ is surjective, since if there is some $e_0 \in E$ such that $e_0 \neq h_{\phi}(f)$ for all $f \in F$, then $h_{\phi}(h_{\psi}(e_0)) \neq e_0$. 
But
\begin{displaymath}
\begin{split}
 a \cdot h_{\phi}(h_{\psi}(e_0)) & = g_{\psi}(g_{\phi}(a)) \cdot e_0 \\
                                 & = a \cdot e_0   
\end{split}
\end{displaymath}
for all $a \in A$.
Thus $e_0$ and $h_{\phi}(h_{\psi}(e_0))$ are distinct duplicate columns, violating $C$'s biextensionality.
So there can be no such $e_0$.
\end{proof}

\bigskip
\noindent \textbf{Claim 39}.
$C$ and $D$ are biextensionally equivalent iff they are homotopy equivalent.

\begin{proof}
We want to show that $C \simeq D$ (interpreted as homotopy equivalence) iff $\hat{C} \cong \hat{D}$.
To do this, we will show that $C \simeq \hat{C}$.
It follows that if $C \simeq D$, then $\hat{C} \simeq C \simeq D \simeq \hat{D}$.
Since we just showed that for biextensional frames, homotopy equivalence and isomorphism coincide, it follows that $\hat{C} \cong \hat{D}$.
Conversely, if $\hat{C} \cong \hat{D}$, then $C \simeq \hat{C} \simeq \hat{D} \simeq D$.

Thus we show that $C \simeq \hat{C}$.
We therefore construct a pair of morphisms $(g, h) : C \to \hat{C}$ and $(j, k) : \hat{C} \to C$.
We define $g : A \to \hat{A}$ by $g(a) = \hat{a}$, and $k : E \to \hat{E}$ by $k(e) = \hat{e}$.
For $h : \hat{E} \to E$ and $j : \hat{A} \to A$, we can send each equivalence class to some member of that class.

We need to show that $(g \circ j, k \circ h)$ is a morphism homotopic to the identity on $\hat{C}$, and likewise that $(j \circ g, h \circ k)$ is a morphism homotopic to the identity on $C$.
Clearly $g \circ j$ and $k \circ h$ are the identity on $\hat{A}$ and $\hat{E}$, respectively, so we are done with the first case.
For the second case, $j \circ g$ and $h \circ k$ need not be the identity, but we see that $j(g(a)) \sim a$ and $h(k(e)) \sim e$ for all $a \in A, e \in E$.
So $j(g(a)) \cdot e = a \cdot e$, which implies that $(j \circ g, \text{id}_E)$ is a morphism.
\end{proof}
\bigskip
\noindent We now prove results from the main text.

\bigskip
\noindent \textbf{Claim 40}.
If $C \cong D$ then $C \simeq D$.

\begin{proof}
If $C \cong D$ then there is some morphism $(g, h) : C \to D$ such that $g$ and $h$ are bijective.
Thus they have inverses which together form a morphism $(g^{-1}, h^{-1})$.
Then these two morphisms compose to the identity in both orders, and $C \simeq D$.
\end{proof}
We now prove Claim 22.

\bigskip
\noindent \textbf{Claim 41}.
The committing definition of additive subagent implies the currying definition.
\begin{proof}
Let $X, Y, Z$ be sets such that $X \subseteq Y$.
Assume $C \simeq (X, Z, \diamond)$ and $D \simeq (Y, Z, \bullet)$,  and let $p: Y \times Z \to W$ be such that $x \diamond z = p(x, z)$ and $y \bullet z = p(y, z)$.

If $D = (B, F, \star)$, then by assumption there are morphisms $(g_0, h_0): D \to (Y, Z, \bullet)$ and $(g_1, h_1) : (Y, Z, \bullet) \to D$ which compose to something homotopic to the identity in both orders.

Now we want to show the existence of some Cartesian frame $M$ such that $C \simeq D^{\circ}(M)$ with $|Env(M)| = 1$.
So we define $M = (X, \{e\}, \cdot)$, where $x \cdot e = g_1(x)$.
Then $D^{\circ}(M) = (X, \{e\} \times F, \star')$ where
\[
\begin{split}
    x \star' (e, f) & = (x \cdot e) \star f \\
                    & = g_1(x) \star f \\
                    & = x  \bullet h_1(f).
\end{split}
\]
We now show that $D^{\circ}(M) \simeq (X, Z, \diamond)$ by constructing morphisms $(g_2, h_2):(X, Z, \diamond) \to D^{\circ}(M)$ and $(g_3, h_3):D^{\circ}(M) \to (X, Z, \diamond)$ which compose to something homotopic to the identity in both orders.
We start simply: let $g_2, g_3$ be the identity on $X$.
Let $h_2(e, f) = h_1(f)$ and $h_3(z) = (e, h_0(z))$.

Next we check that these are indeed morphisms:
\[
\begin{split}
    g_2(x) \star' (e, f) & = x \star' (e, f) \\
                         & = x \bullet h_1(f) \\
                         & = x \diamond h_1(f) \\
                         & = x \diamond h_2(e, f).
\end{split}
\]
Also,
\[
\begin{split}
    g_3(x) \diamond z & = x \diamond z \\
                      & = x \bullet z \\
                      & = x \bullet h_1(h_0(z)) \\
                      & = x \star' (e, h_0(z)) \\
                      & = x \star' h_3(z).
\end{split}
\]
We also see that the two morphisms compose to something homotopic to the identity in both orders, since $g_2 \circ g_3 = \text{id}_X = g_3 \circ g_2$.

Thus we see $C \simeq (X, Z, \diamond) \simeq D^{\circ}(M)$, and $|Env(M)| = 1$, concluding the proof.
\end{proof}

\bigskip
\noindent \textbf{Claim 42}. The currying definition of additive subagent implies the committing definition.

\begin{proof}
Assume $C \simeq D^{\circ}(M)$ with $|Env(M)| = 1$.
Let $D = (Y, Z, \bullet)$, and let $f: Y \times Z \to W$ be given by $f(y, z) = y \bullet z$.

Let $X \subseteq Y$ be given by $X = Image(M)$.
Since $|Env(M)| = 1$, we have $M \simeq \bot_X$.
So, $C \simeq D^{\circ}(M) \simeq D^{\circ}(\bot_X)$.

But $D^{\circ}(\bot_X) = (X, \{j\} \times Z, \cdot)$ where $x \cdot (j, z) = f(x, z)$, and hence is isomorphic to $(X, Z, \diamond)$ where $x \diamond z = f(x, z)$.
So $C \simeq D^{\circ}(M) \simeq D^{\circ}(\bot_X) \cong (X, Z, \diamond)$;
that is, $C \simeq (X, Z, \diamond)$ and $D = (Y, Z, \bullet)$ with $X \subseteq Y$, as in the committing definition.
\end{proof}

\bigskip
\noindent \textbf{Claim 43}.
The externalizing definition of multiplicative subagent implies the currying definition.

\begin{proof}
Assume three sets $X, Y, Z$, and a function $p:X \times Y \times Z \to W$ such that $C \simeq (X, Y \times Z, \diamond)$ and $D \simeq (X \times Y, Z, \bullet)$ where $x \diamond (y, z) = (x, y) \bullet z = p(x, y, z)$.

Let $D = (B, F, \star)$.
Then by assumption there are morphisms $(g_0, h_0) : D \to (X \times Y, Z, \bullet)$ and $(g_1, h_1) : (X \times Y, Z, \bullet) \to D$ which compose to something homotopic to the identity in both orders.

We now define a Cartesian frame $M$ over $B$ by $M = (X, Y \times B, \cdot)$ where, if $g_0(b)=(x, y)$, then $x \cdot (y, b) = b$, and otherwise $x \cdot (y, b) = g_1(x, y)$.
Clearly $Image(M) \subseteq B$;
to show equality, consider that for any $b \in B$, if we let $g_0(b) = (x, y)$, then $x \cdot (y, b) = b$.

Clearly $(x \cdot (y, b)) \star f = g_1(x, y) \star f$ when $(x, y) \neq g_0(b)$.
As it turns out, this equality holds even when $(x, y) = g_0(b)$:
\[
\begin{split}
    (x \cdot (y, b)) \star f & = b \star f \\
                             & = g_1(g_0(b)) \star f \\
                             & = g_1(x, y) \star f.
\end{split}
\]
So now we have $D^{\circ}(M) = (X, Y \times B \times F, \star')$ where
\[
\begin{split}
    x \star' (y, b, f) & = (x \cdot (y, b)) \star f \\
                       & = g_1(x, y) \star f \\
                       & = (x, y) \bullet h_1(f).
\end{split}
\]
Our final step is to show that $(X, Y \times Z, \diamond) \simeq D^{\circ}(M)$.
To do so we construct a pair of morphisms $(g_2, h_2) : (X, Y \times Z, \diamond) \to D^{\circ}(M)$ and $(g_3, h_3) : D^{\circ}(M) \to (X, Y \times Z, \diamond)$ that compose to something homotopic to the identity in both orders.

We can let $g_2, g_3 = \text{id}_X$.
Let $h_2(y, b, f) = (y, h_1(f))$, and let $h_3(y, z) = (y, b, h_0(z))$ for some $b \in B$ (it does not matter which).
We check that these constructions are morphisms:
\[
\begin{split}
    g_2(x) \star' (y, b, f) & = x \star' (y, b, f) \\
                            & = (x, y) \bullet h_1(f) \\
                            & = p(x, y, h_1(f)) \\
                            & = x \diamond (y, h_1(f)) \\
                            & = x \diamond h_2(y, b, f). 
\end{split}
\]
And the other pair:
\[
\begin{split}
    g_3(x) \diamond (y, z) & = x \diamond (y, z) \\
                           & = p(x, y, z) \\
                           & = (x, y) \bullet z \\
                           & = (x, y) \bullet h_1(h_0(z)) \\
                           & = x \star' (y, b, h_0(z)) \\
                           & = x \star' h_3(y, z).
\end{split}
\]
Finally, $(g_2, h_2)$ and $(g_3, h_3)$ compose to something homotopic to the identity in both orders since $g_2, g_3$ are the identity on $X$.
Thus, $C \simeq (X, Y \times Z, \diamond) \simeq D^{\circ}(M)$ and $Image(M) = Agent(D)$, satisfying the currying definition.
\end{proof}

\bigskip
\noindent \textbf{Claim 44}.
The currying definition of multiplicative subagent implies the externalizing definition.

\begin{proof}
Assume $C \simeq D^{\circ}(M)$, where $M$ is a Cartesian frame over $Agent(D)$ with $Image(M) = Agent(D)$.
Explicitly, let $M = (X, Y, \cdot)$ and $D = (B, Z, \star)$, and let $f : X \times Y \times Z \to W$ be given by $f(x, y, z) = (x \cdot y) \star z$.

Thus $D^{\circ}(M) = (X, Y \times Z, \diamond)$ where $x \diamond (y, z) = (x \cdot y) \star z = f(x, y, z)$.
It suffices to show that $D \simeq (X \times Y, Z, \bullet)$ where $(x, y) \bullet z = f(x, y, z)$.
To do so we construct morphisms $(g_0, h_0) : D \to (X \times Y, Z, \bullet)$ and $(g_1, h_1) : (X \times Y, Z, \bullet) \to D$ that compose to something homotopic to the identity in both orders.

We can let $h_0, h_1$ be the identity on $Z$.
Let $g_1(x, y) = x \cdot y$.
We see that $g_1$ is surjective since $Image(M) = B$, so $g_1$ has a right inverse.
Let $g_0$ be any right inverse to $g_1$, so $g_1(g_0(b))=b$ for all $b \in B$.
We check that these are indeed morphisms:
\[
\begin{split}
    g_1(x, y) \star z & = (x \cdot y) \star z \\
                      & = f(x, y, z) \\
                      & = (x, y) \bullet z \\
                      & = (x, y) \bullet h_1(z).
\end{split}
\]
And also, if we let $g_0(b) = (x, y)$:
\[
\begin{split}
    g_0(b) \bullet z & = (x, y) \bullet z \\
                     & = f(x, y, z) \\
                     & = (x \cdot y) \star z \\
                     & = g_1(x \cdot y) \star z \\
                     & = g_1(g_0(b)) \star z \\
                     & = b \star h_0(z).
\end{split}
\]
We see that $(g_0, h_0)$ and $(g_1, h_1)$ compose to something homotopic to the identity since $h_0, h_1$ are the identity on $Z$.
Thus $D \simeq (X \times Y, Z, \bullet)$, completing the proof.
\end{proof}

\noindent \textbf{Claim 45}. Let $C=(A, E, \cdot)$. Then
\begin{enumerate}[label=(\arabic*)]
    \item for all partitions $B$ of $A$, $External^B (C) \triangleleft_{\times} C$ and $External^{/B} (C) \triangleleft_{\times} C$.
    \item for all partitions $F$ of $E$, $C \triangleleft_{\times} Internal^F (C)$ and $C \triangleleft_{\times} Internal^{/F} (C)$.
\end{enumerate}

\begin{proof}
(1) As one might expect, this result is nearly immediate from the externalizing definition of multiplicative subagent.
Let $External^B (C) = (A/B, B \times E, \star)$ where $q \star (b, e) = q(b) \cdot e$.
Consider the Cartesian frame $(A/B \times B, E, \diamond)$ where $(q, b) \diamond e = q(b) \cdot e$.
Clearly then we have our three sets, $A/B, B, E$ with a function $f : A/B \times B \times E \to W$ given by $f(q, b, e) = (q, b) \diamond e = q(b) \cdot e = q \star (b, e)$.
So we need that $C \simeq (A/B \times B, E, \diamond)$.

So we construct two morphisms $(g_0, h_0) : C \to (A/B \times B, E, \diamond)$ and $(g_1, h_1): (A/B \times B, E, \diamond) \to C$ which compose to something homotopic to the identity in both orders.
We can let $h_0, h_1$ be the identity on $E$.
Let $g_0 : A \to A/B \times B$ be given by $g_0(a) = (q, b)$ where $q(b) = a$.
Let $g_1: A/B \times B \to A$ be given by $g_1(q, b) = q(b)$.
Then we have
\[
\begin{split}
    g_0(a) \diamond e & = (q, b) \diamond e \\
                      & = q(b) \cdot e \\
                      & = a \cdot h_0(e)
\end{split}
\]
and
\[
\begin{split}
    g_1(q, b) \cdot e & = q(b) \cdot e \\
                      & = (q, b) \diamond e \\
                      & = (q, b) \diamond h_1(e).
\end{split}
\]
Clearly the two morphisms compose to something homotopic to the identity in both orders since $h_0, h_1$ are the identity on $E$.\footnote{In fact these morphisms are bijective and so establish an isomorphism, as the reader can verify.}
The proof that $External^{/B}(C) \triangleleft_{\times} C$ is similar.

\bigskip
(2) This follows from the fact that $(Internal^F(C))^* \cong External^F(C^*) \triangleleft_{\times} C^*$ and $(Internal^{/F}(C))^* \cong External^{/F}(C^*) \triangleleft_{\times} C^*$ and the fact that multiplicative subagent is equivalent to multiplicative sub-environment.
\end{proof}

\section{Categorical Framework}

This appendix develops the category-theoretic perspective on Cartesian frames.
Here we show how to develop the CF constructions found in the main body via categorical definitions.
We also prove the equivalence of the category-theoretic constructions and those found in the main body.

As was stated above, $\text{Chu}(W)$ is the category of Chu spaces (Cartesian frames) defined over the set $W$, having as arrows all morphisms between Chu spaces, with composition defined as in the main text.
To prove that $\text{Chu}(W)$ is indeed a category, we need only show that composition of morphisms is well-defined and associative and that there exist identity morphisms.

\begin{proof}
For identity morphisms, $(\text{id}_A, \text{id}_E)$ is clearly an identity on $C = (A, E, \cdot)$, where $\text{id}_X$ is the identity map from $X$ to itself.

As defined in the main body, given Cartesian frames $C_i = (A_i, E_i, \cdot_i)$ and morphisms $(g_0, h_0) : C_0 \to C_1$ and $(g_1, h_1) : C_1 \to C_2$, their composition, $(g_1, h_1) \circ (g_0, h_0)$ is given by $(g_1 \circ g_0, h_0 \circ h_1) : C_0 \to C_2$.
Then, for all $a_0 \in A_0$ and all $e_2 \in E_2$:
\begin{displaymath}
\begin{split}
  a_0 \cdot_0 h_0(h_1(e_2)) & = g_0(a_0) \cdot_1 h_1(e_2) \\
                            & = g_1(g_0(a_0)) \cdot_2 e_2
\end{split}
\end{displaymath}
Thus $g_1 \circ g_0$ and $h_0 \circ h_1$ satisfy the adjointness condition, and $(g_1 \circ g_0, h_0 \circ h_1)$ is a morphism.

Associativity follows from the fact that the new morphism is a pair of compositions of set functions, and composition is associative for set functions.
\end{proof}

\noindent So $\text{Chu}(W)$ is indeed a category.
In fact, $\text{Chu}(W)$ is a \emph{self-dual} category; that is, $\text{Chu}(W)$ is isomorphic to its dual.\footnote{Chu$(W)$ is known as a \emph{*-autonomous} category, systematically studied by \cite{barr1979autonomous}.}

\bigskip
\noindent \textbf{Claim 46}.
$-^*$ is an isomorphism between $\text{Chu}(W)$ and $\text{Chu}(W)^{\text{op}}$.

\begin{proof}
First we must show that $-^*$ is a functor.
By definition, the objects of $\text{Chu}(W)$ are the same as those of $\text{Chu}(W)^{\text{op}}$, and the arrows from $D$ to $C$ in $\text{Chu}(W)^{\text{op}}$ are the arrows from $C$ to $D$ in $\text{Chu}(W)$, with composition in reverse order.
We see that $-^*$ preserves identity arrows.

We now show that $-^*$ preserves composition.
To this end, we have
\begin{displaymath}
\begin{split}
    (g_0, h_0)^* \circ^{\text{op}} (g_1, h_1)^* & = (h_1, g_1) \circ (h_0, g_0) \\
                                                & = (h_1 \circ h_0, g_0 \circ g_1) \\
                                                & = ((g_0, h_0) \circ (g_1, h_1))^*
\end{split}
\end{displaymath}
We now show that it is an isomorphism.
If we abuse notation and also write $-^* : \text{Chu}(W)^{\text{op}} \to \text{Chu}(W)$ as the functor given by $(E, A, \star)^* = (A, E, \cdot)$, where $a \cdot e = e \star a$, then clearly $-^* : \text{Chu}(W) \to \text{Chu}(W)^{\text{op}}$ and $-^* : \text{Chu}(W)^{\text{op}} \to \text{Chu}(W)$ compose to the identity in both orders.
\end{proof}

\noindent For our purposes, the claim that $\text{Chu}(W)$ is self-dual really just amounts to the claim that the transpose of a Cartesian frame is itself a Cartesian frame.
As a particular case, we have $(C^*)^* = C$.

$\text{Chu}(W)$ also has initial and terminal objects:

\bigskip
\noindent \textbf{Definition 47}.
Let $0 = (\{\}, \{e\}, \cdot)$, where $\text{Agent}(0)$ is empty, $\{e\}$ is any singleton, and $\cdot$ is trivial.
Let $\top = (\{a\}, \{\}, \cdot)$, where $\text{Env}(0)$ is empty, $\{a\}$ is any singleton, and $\cdot$ is trivial.

\bigskip
\noindent \textbf{Claim 48}.
$0$ is initial in $\text{Chu}(W)$ and $\top$ is terminal in $\text{Chu}(W)$.

\begin{proof}
Let $C = (A, E, \cdot)$ be a Cartesian frame over $W$.
Then a morphism from $0$ to $C$ is a pair of functions $(g, h)$, where $g$ is a function from $\{\}$ to $A$ and $h$ is a function from $E$ to $\{e\}$.
Considered as objects in \textbf{Set}, the category of sets, $\{\}$ is initial, and any singleton is terminal.
So $g$ and $h$ are unique, and $0$ is initial.
The proof that $\top$ is terminal is nearly identical.
\end{proof}

\noindent Technically these terms are ill-defined since the denoted frame is not unique (there exists one for each singleton).
The obvious solution is to show that each initial and terminal object is isomorphic (a trivial corollary), take the equivalence classes of initial and terminal objects up to isomorphism, and let $0$ and $\top$ denote a member of these classes.

In addition to our definition of $\bot_S$, we require one more Cartesian frame for what follows:

\bigskip
\noindent \textbf{Definition 49}.
Given $S \subseteq W$, $1_S$ is the Cartesian frame $1_S = (\{a\}, S, \star)$, where $\star$ is given by $a \star s = s$ for all $s \in S$.
$1$ is the Cartesian frame $1_W$.

\bigskip
\noindent Note that $(\bot_S)^* = 1_S$ and $(1_S)^* = \bot_S$ for any $S \subseteq W$.

We introduce a final definition of subagent for ease of use in proofs.
It will be shown to be equivalent to the other two.

\bigskip
\noindent \textbf{Definition 50} (Covering definition of subagent).
Let $C = (A, E, \cdot), D = (B, F, \star)$ be Cartesian frames over $W$.
We say that $C \triangleleft D$ if for all $e \in E$, there is a morphism $(g, h) : C \to D$ and an $f \in F$ such that $e = h(f)$.

\noindent We call this the covering definition because $E$ is covered by $F$, or equivalently, there is a morphism $(g, h)$ with $h$ surjective.

We're now in a position to prove that the categorical and currying definitions are equivalent.
We'll do this by a series of implications:
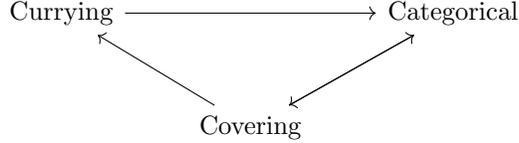
\[\begin{tikzcd}[row sep=large]
	{\text{Currying}} && {\text{Categorical}} \\
	& {\text{Covering}}
	\arrow[from=2-2, to=1-1]
	\arrow[from=1-1, to=1-3]
	\arrow[from=1-3, to=2-2]
	\arrow[from=2-2, to=1-3]
\end{tikzcd}\]

\bigskip
\noindent \textbf{Claim 51}.
The covering and categorical definitions of subagent are equivalent.

\begin{proof}
Recall that each arrow $C \to \bot$ corresponds uniquely to an element $e$ of the environment $E$ of that frame; we denote the morphism associated with $e$ by $\phi_e$.
If $C \triangleleft D$ according to the categorical definition, then, for every $\phi_e: C \to \bot$, there are arrows $(g, h) : C \to D$, $\psi_f : D \to \bot$ such that $\phi_e = \psi_f \circ (g, h)$.
But $\psi_f \circ (g, h)$ sends $j$ to $h(\psi_f(j)) = h(f)$, and thus equals $e$ iff $h(f) = e$.
This is just the covering definition, so the two definitions are equivalent.
\end{proof}

\bigskip
\noindent \textbf{Claim 52}.
The covering definition of subagent implies the currying definition.

\begin{proof}
Let $C=(A, E, \cdot)$ and $D=(B, F, \star)$ be Cartesian frames such that $C \triangleleft D$ according to the covering definition.

Let $X = \text{hom}(C, D)$, and let $Z=(A, X, \diamond)$ be a Cartesian frame over $B$, with $\diamond$ given by $a \diamond (g, h) = g(a)$.
Thus $D^{\circ}(Z)=(A, X \times F, \bullet)$, where 
\begin{align*}
    a \bullet ((g, h), f) & = (a \diamond (g, h)) \star f \\
                          & = g(a) \star f
\end{align*}
for all $a \in A$, $(g, h) \in X$, and $f \in F$.

We now want to show that $C \simeq D^{\circ}(Z)$.
To do so, we will construct morphisms $(g_0, h_0): C \to D^{\circ}(Z)$ and $(g_1, h_1): D^{\circ}(Z) \to C$ which compose to something homotopic to the identity in both orders.

We can let $g_0, g_1$ be the identity on $A$.
We let $h_0: X \times F \to E$ be defined by $h_0((g,h), f) = h(f)$.
Finally, we let $h_1 : E \to X \times F$ be given by $h_1(e) = ((g, h), f)$ such that $h(f)=e$.
By the covering definition, there is always such a morphism $(g, h)$ and an $f \in F$.

We now check that $(g_0, h_0)$ is a morphism.
We have
\begin{align*}
    g_0(a) \bullet ((g, h), f) & = a \bullet ((g, h), f) \\
                               & = g(a) \star f \\
                               & = a \cdot f \\
                               & = a \cdot h_0((g, h), f). \\
\end{align*}
Similarly, for $(g_1, h_1)$ we have
\begin{align*}
    g_1(a) \cdot e & = a \cdot e \\
                   & = a \cdot h(f) \\
                   & = g(a) \star f \\
                   & = a \bullet ((g, h), f) \\
                   & = a \bullet h_1(e)
\end{align*}
since $h_1(e) = ((g, h), f)$, where $h(f)=e$. Clearly $(g_0, h_0)$ and $(g_1, h_1)$ compose to something homotopic to the identity in both orders, since $g_0, g_1$ are the identity on $A$.

\end{proof}

\bigskip
\noindent \textbf{Claim 53}.
The currying definition of subagent implies the categorical definition.

\begin{proof}
Let $C = (A, E, \cdot)$, $D = (B, F, \star)$ be Cartesian frames over $W$.
Let $Z = (X, Y, \diamond)$ be a Cartesian frame over $B$ such that $C \simeq D^{\circ}(Z)$.
Explicitly, $D^{\circ}(Z) = (X, Y \times F, \bullet)$ where $x \bullet (y, f) = (x \diamond y) \star f$.
Then there are morphisms $\sigma : C \to D^{\circ}(Z)$ and $\tau : D^{\circ}(Z) \to C$ which compose to something homotopic to the identity in both orders.

Let $\phi_e : C \to \bot$ be a morphism corresponding to $e \in E$.
We want to show that there are morphisms $\xi : C \to D$, $\psi : D \to \bot$ such that $\phi_e = \psi \circ \xi$.
First, note that for any $y \in Y$ we can define a morphism $(g_y, h_y)$ from $D^{\circ}(Z)$ to $D$ by
\[g_y(x) = x \diamond y\]
\[h_y(f) = (y, f).\]
We check that $(g_y, h_y)$ is indeed a morphism:
\[
\begin{split}
    g_y(x) \star f & = (x \diamond y) \star f \\
                   & = x \bullet (y, f) \\
                   & = x \bullet h_y(f).
\end{split}
\]
Thus $(g_y, h_y) \circ \sigma$ is a morphism from $C$ to $D$ for any $y \in Y$.

We now construct a $\psi : D \to \bot$ such that $\psi \circ ((g_y, h_y) \circ \sigma) = \phi_e$, i.e., such that the following diagram commutes:
\begin{center}
    \begin{tikzcd}[column sep=large]
C \arrow[r, "\phi_e"] \arrow[d, "\sigma"] & \bot                \\
D^{\circ}(Z) \arrow[r, "{(g_y, h_y)}"]  & D \arrow[u, "\psi", swap]
    \end{tikzcd}
\end{center}
This will show that $\phi_e$ factors through $D$, establishing our result.
To do so we will make use of the fact that $\sigma$ and $\tau$ compose to something homotopic to the identity in both orders.
If we let $\sigma = (g_\sigma, h_\sigma )$ and $\tau = (g_\tau, h_\tau )$ then we know that $h_\tau (h_{\phi_e} (j)) = (y, f)$ for some $(y, f) \in Y \times F$.
Holding this $(y, f)$ fixed, we may select the corresponding morphisms $(g_y, h_y)$ and $\psi_f$.
Then $h_y (h_{\psi_f} (j)) = (y, f)$.
Thus
\[
\begin{split}
  g_{\phi_e}(g_\tau(x)) \cdot_\bot j & = x \bullet h_\tau(h_{\phi_e}(j)) \\
                                     & = x \bullet (y, f) \\
                                     & = x \bullet h_y(h_{\psi_f}(j)) \\
                                     & = g_{\psi_f}(g_y(x)) \cdot_\bot j
\end{split}
\]
so $\phi_e \circ \tau = \psi_f \circ (g_y, h_y)$. This implies that $\phi_e \circ \tau \circ \sigma = \psi_f \circ (g_y, h_y) \circ \sigma$.
Finally, we show that $\phi_e = \phi_e \circ \tau \circ \sigma$:
\[
\begin{split}
  g_{\phi_e} (g_{\tau} (g_{\sigma} (a))) \cdot_\bot j & = a \cdot (h_{\sigma} (h_\tau (h_{\phi_e} (j)))) \\
                                                     & = a \cdot e \\
                                                     & = g_{\phi_e}(a) \cdot_\bot j. 
\end{split}
\]
Thus $\phi_e  = \psi_f \circ (g_y, h_y) \circ \sigma$, and so $\phi_e$ factors through $D$.
\end{proof}

\noindent So we see that the three definitions are equivalent.
We can also give categorical definitions of additive and multiplicative subagents, which are, as one might expect, quite similar to the above categorical definition of subagent:

\bigskip
\noindent \textbf{Definition 54} (Categorical definition of additive subagent).
We say $C \triangleleft_+ D$ if there exists a single morphism $\phi_0 : C \to D$ such that for every morphism $\phi : C \to \bot$ there exists a morphism $\phi_1: D \to \bot$ such that $\phi$ is homotopic to $\phi_1 \circ \phi_0$.

\bigskip
\noindent Note that the difference between this definition and the categorical definition of (general) subagent is that, in the case of \emph{additive} subagents, we require the morphism from $C$ to $D$ to be unique, and we require only \emph{homotopy}, not equality.

\bigskip
\noindent \textbf{Claim 55}.
The categorical definition of additive subagent implies the committing definition.

\begin{proof}
Let $C = (A, E, \cdot)$ and $D = (B, F, \bullet)$ be Cartesian frames over $W$.
Let $(g_0, h_0) : C \to D$ be such that for all $(g, h): C \to \bot$, there exists a $(g', h'):D \to \bot$ such that $(g', h') \circ (g_0, h_0)$ is homotopic to $(g, h)$.
As usual, let $\bot = (W, \{j\}, \star)$.

Let $X = \{g_0(a) : a \in A\} \subseteq B$.
Let $p : B \times F \to W$ be given by $p(b, f) = b \bullet f$.
Then we must show that $C \simeq (X, F, \diamond)$ where $x \diamond z = p(x, f)$ (i.e., $(X, F, \diamond)$ is just $D$ with a restricted agent).
To show homotopy equivalence we must construct $(g_1, h_1) : C \to (X, F, \diamond)$ and $(g_2, h_2): (X, F \diamond) \to C$ which compose to something homotopic to the identity in both orders.

So we define $g_1 : A \to X$ by $g_1(a) = g_0(a)$.
$g_1$ is surjective and thus has a right inverse.
We let $g_2: X \to A$ be any right inverse to $g_1$.
We let $h_1 : F \to E$ be $h_1(z) = h_0(z)$.
So $(g_1, h_1)$ is $(g_0, h_0)$ but with codomain $(X, F, \diamond)$.

We define $h_2: E \to F$ by way of the fact that morphisms into $\bot$ are completely determined by elements of the environment of the domain.
Thus given an $e \in E$, let $\phi_e = (g_e, h_e) : C \to \bot$ be the corresponding morphism.
Then by assumption there is a morphism $(g'_e, h'_e): D \to \bot$ such that $(g'_e, h'_e) \circ (g_0, h_0)$ is homotopic to $(g_e, h_e)$.
We define $h_2$ by $h_2(e) = h'_e (j)$.

Now we show that these do indeed form morphisms.
$(g_1, h_1)$ is a morphism since $(g_0, h_0)$ is.
For $(g_2, h_2)$ we see that
\[
\begin{split}
    x \diamond h_2(e) & = g_1(g_2(x)) \diamond h'_e(j) \\
                      & = g'_e(g_0(g_2(x))) \star j \\
                      & = g_2(x) \cdot h_e(j) \\
                      & = g_2(x) \cdot e.
\end{split}
\]
As our final step we check that these two morphisms compose to something homotopic to the identity in both orders.
Clearly $(g_1, h_1) \circ (g_2, h_2)$ does since $g_1 \circ g_2 = \text{id}_X$.
In the other order we have
\[
\begin{split}
    g_2(g_1(a)) \cdot e & = g_1(a) \diamond h_2(e) \\
                        & = g_0(a) \diamond h'_e(j) \\
                        & = g'_e(g_0(a)) \star j \\
                        & = a \cdot h_e(j) \\
                        & = a \cdot e.
\end{split}
\]
So $C \simeq (X, F, \diamond)$, and $C \triangleleft_{+} D$ according to the committing definition.
\end{proof}

\bigskip
\noindent \textbf{Claim 56}. The committing definition of additive subagent implies the categorical definition.

\begin{proof}
Let $X, Y$, and $Z$ be sets such that $X \subseteq Y$, let $p: Y \times Z \to W$, and let $C \simeq (X, Z, \diamond)$ and $D \simeq (Y, Z, \bullet)$.

By assumption there are morphisms $(g_1, h_1): C \to (X, Z, \diamond)$ and $(g_2, h_2) : (X, Z, \diamond) \to C$ which compose to something homotopic to the identity in both orders, and similarly there are morphisms $(g_3, h_3) : D \to (Y, Z, \bullet)$, $(g_4, h_4) : (Y, Z, \bullet) \to D$ which compose to something homotopic to the identity in both orders.
Further, we define $(g_0, h_0) : (X, Z, \diamond) \to (Y, Z, \bullet)$ by letting $g_0$ be the inclusion $X \xhookrightarrow{} Y$ and $h_0 = \text{id}_Z$.

We define our single morphism $\phi : C \to D$ by $\phi = (g_4, h_4) \circ (g_0, h_0) \circ (g_1, h_1)$, so the following diagram commutes:
\[
\begin{tikzcd}
C \arrow[r, "\phi"] \arrow[d, "{(g_1, h_1)}"'] & D                                            \\
{(X, Z, \diamond)} \arrow[r, "{(g_0, h_0)}"']  & {(Y, Z, \bullet)} \arrow[u, "{(g_4, h_4)}"']
\end{tikzcd}
\]
Our goal now is to show that for any morphism $(g, h): C \to \bot$, we can construct a morphism $(g', h') : D \to \bot$ such that $(g, h)$ is homotopic to $(g', h') \circ \phi$.

So let $\bot = (W, \{j\}, \star)$, $C = (A, E, \cdot_0)$, $D = (B, F, \cdot_1)$.
We define $h': \{j\} \to F$ by $h' = h_3 \circ h_2 \circ h$.
Let $g': B \to W$ be given by $g'(b) = b \cdot_1 h'(j)$.
Then $(g', h')$ is a morphism:
\[
\begin{split}
    g'(b) \star j & = g'(b) \\
                  & = b \cdot_1 h'(j). 
\end{split}
\]
Our final step is to show that $(g, h)$ is indeed homotopic to $(g', h') \circ \phi = (g', h') \circ (g_4, h_4) \circ (g_0, h_0) \circ (g_1, h_1)$.
It suffices to show that $(g, h_1 \circ h_0 \circ h_4 \circ h')$ is a morphism;
equivalently, since $h_0$ is the identity and $h' = h_3 \circ h_2 \circ h$, we show that
\[
(g, h_1 \circ h_4 \circ h_3 \circ h_2 \circ h)
\]
is a morphism:
\[
\begin{split}
    g(a) \star j & = a \cdot_0 h(j) \\
                 & = a \cdot_0 h_1(h_2(h(j))) \\
                 & = g_1(a) \diamond h_2(h(j)) \\
                 & = g_1(a) \bullet h_2(h(j)) \\
                 & = g_1(a) \bullet h_4(h_3(h_2(h(j)))) \\
                 & = g_1(a) \diamond h_4(h_3(h_2(h(j)))) \\
                 & = a \cdot_0 h_1(h_4(h_3(h_2(h(j))))). \\
\end{split}
\]
\end{proof}

\noindent \textbf{Definition 57} (Categorical definition of multiplicative subagent).
We say $C \triangleleft_{\times} D$ if for every morphism $\phi: C \to \bot$, there exist morphisms $\phi_0 : C \to D$, $\phi_1 : D \to \bot$ such that $\phi = \phi_1 \circ \phi_0$, and for every morphism $\psi : 1 \to D$, there exist morphisms $\psi_0 : 1 \to C$, $\psi_1 : C \to D$ such that $\psi = \psi_1 \circ \psi_0$.
That is, the following diagrams commute:
\[
 \begin{tikzcd}[row sep=large, column sep=large]
C \arrow[r, "\phi"] \arrow[d, "\phi_0"] & \bot \\ 
D \arrow[ru, "\phi_1", swap]
\end{tikzcd}
\begin{tikzcd}[row sep=large, column sep=large]
                                          & C \arrow[d, "\psi_1"] \\
1 \arrow[r, "\psi"'] \arrow[ru, "\psi_0"] & D                    
\end{tikzcd}
\]

\bigskip
\noindent We will show that this definition is equivalent to the others by way of yet another equivalent definition, this one in terms of sub-environments:

\bigskip
\noindent \textbf{Definition 58} (Sub-environment definition of multiplicative subagent).
We say $C \triangleleft_{\times} D$ if $C \triangleleft D$ and $C \triangleleft^* D$. Equivalently, $C \triangleleft_{\times} D$ if $C \triangleleft D$ and $D^* \triangleleft C^*$.

\bigskip
\noindent \textbf{Claim 59}.
The categorical definition of multiplicative subagent is equivalent to the sub-environment definition.

\begin{proof}
Clearly both definitions require that $C \triangleleft D$.

The condition that for every morphism $\psi: 1 \to D$, there are morphisms $\psi_0 : 1 \to C$ and $\psi_1 : C \to D$ such that $\psi = \psi_1 \circ \psi_0$ is equivalent to the condition that for every morphism $\psi^* : D^* \to \bot$ there are morphisms $\psi^*_0 : C^* \to \bot$ and $\psi^*_1 : D^* \to C^*$ such that $\psi^* = \psi^*_0 \circ \psi^*_1$.
This is the categorical definition of $D^* \triangleleft C^*$.
\end{proof}

\bigskip
\noindent \textbf{Claim 60}.
The sub-environment definition of multiplicative subagent is equivalent to the externalizing definition.

\begin{proof}
First we show that the sub-environment definition implies the externalizing definition.
Let $C = (A, E, \cdot)$ and $D = (B, F, \star)$ be Cartesian frames over $W$ with $C \triangleleft D$ and $C \triangleleft^* D$.
We define a function $p: A \times \text{hom}(C, D) \times F \to W$ by
\[
\begin{split}
    p(a, (g, h), f) & = g(a) \star f \\
                    & = a \cdot h(f).
\end{split}
\]
We want to show that $C \simeq (A, \text{hom}(C, D) \times F, \diamond)$ and $D \simeq (A \times \text{hom}(C, D), F, \bullet)$, where $a \diamond ((g, h), f) = p(a, (g, h), f) = (a, (g, h)) \bullet f$.

So we construct morphisms $(g_0, h_0) : C \to (A, \text{hom}(C, D) \times F, \diamond)$ and $(g_1, h_1) : (A, \text{hom}(C, D) \times F, \diamond) \to C$ which compose to something homotopic to the identity in both orders.
Let $g_0, g_1$ be the identity on $A$.
Let $h_0$ be defined by $h_0((g, h), f) = h(f)$.
By the covering definition of subagent, $h_0$ is surjective, and has a right inverse.
Let $h_1$ be any right inverse to $h_0$.

We check that these are indeed morphisms:
\[
\begin{split}
    g_0(a) \diamond ((g, h), f) & = a \diamond ((g, h), f) \\
                                & = p(a, (g, h), f) \\
                                & = a \cdot h(f) \\
                                & = a \cdot h_0((g, h), f)
\end{split}
\]
and, if $h_1(e) = ((g, h), f)$: 
\[
\begin{split}
    g_1(a) \cdot e & = a \cdot h_0(h_1(e)) \\
                   & = a \cdot h_0((g, h), f) \\
                   & = a \cdot h(f) \\
                   & = p(a, (g, h), f) \\
                   & = a \diamond ((g, h), f) \\
                   & = a \diamond h_1(e).
\end{split}
\]
We also see that the two morphisms compose to something homotopic to the identity in both orders since $g_0, g_1$ are the identity on $A$.
So $C \simeq (A, \text{hom}(C, D) \times F, \diamond)$.

We now want to show that $D \simeq (A \times \text{hom}(C, D), F, \bullet)$.
To do so we construct two morphisms $(g_2, h_2) : D \to (A \times \text{hom}(C, D), F, \bullet)$ and $(g_3, h_3) : (A \times \text{hom}(C, D), F, \bullet) \to D$ which compose to something homotopic to the identity in both orders.

We can let $h_2$ and $h_3$ be the identity on $F$.
We define $g_3$ by $g_3(a, (g, h)) = g(a)$.
By the covering definition of subagent and the fact that $D^* \triangleleft C^*$, $g_3$ is surjective and so has a right inverse.
Let $g_2$ be some right inverse to $g_3$.

We check that these are indeed morphisms:
\[
\begin{split}
    g_3(a, (g, h)) \star f & = g(a) \star f \\
                           & = p(a, (g, h), f) \\
                           & = (a, (g, h)) \bullet f \\
                           & = (a, (g, h)) \bullet h_3(f)
\end{split}
\]
and, if $g_2(b) = (a, (g, h))$:
\[
\begin{split}
    g_2(b) \bullet f & = (a, (g, h)) \bullet f \\
                     & = p(a, (g, h), f) \\
                     & = g(a) \star f \\
                     & = g_3(a, (g, h)) \star f \\
                     & = g_3(g_2(b)) \star f \\
                     & = b \star h_2(f).
\end{split}
\]
Clearly the two morphisms compose to something homotopic to the identity in both orders since $h_2, h_3$ are the identity on $F$.
Thus $D \simeq (A \times \text{hom}(C, D), F, \bullet)$, so $C$ and $D$ satisfy the externalizing definition of multiplicative subagent.

Conversely, if $C \triangleleft_{\times} D$ according to the externalizing definition, then we also have $D^* \triangleleft_{\times} C^*$.
However, multiplicative subagent is a stronger notion than subagent, so we have $C \triangleleft D$ and $D^* \triangleleft C^*$, or equivalently $C \triangleleft^* D$.
\end{proof}

\noindent In the main body we chose to present Cartesian frames from a primarily matrix-centered point of view.
We made this decision because matrices are more familiar to the working decision theorist.
But the category-theoretic constructions may offer powerful insights into the underlying mathematics.
For this reason we are interested in rewriting the CF framework from a purely categorical point of view.
We hope to pursue this project in the future.

\printbibliography

\end{document}